\documentclass{article}
\usepackage[T2A]{fontenc}
\usepackage[cp1251]{inputenc}
\usepackage[english,russian]{babel}
\usepackage[tbtags]{amsmath}
\usepackage{amsfonts,amssymb,mathrsfs,amscd}
\usepackage[hyper]{msb-a}
\JournalName{МАТЕМАТИЧЕСКИЙ СБОРНИК}
\JournalName{}
\numberwithin{equation}{section}
\theoremstyle{plain}

\theoremstyle{definition}

\begin{document}

\title{Parabolicity of degenerate singularities of axisymmetric Zhukovsky case}

\author[V.\,A. Kibkalo]{V.\,A.~Kibkalo}
\address{Lomonosov Moscow State University, Moscow, Russia; Moscow Center for Fundamental and Applied Mathematics}
\email{slava.kibkalo@gmail.com}

\udk{517.938.5}

\maketitle

\begin{fulltext}

\begin{abstract}
The degenerate singularities of systems from one well-known multiparameter family of integrable systems of rigid body dynamics are studied. Axisymmetric Zhukovsky systems are considered, i.e. axisymmetric Euler tops after adding a constant gyrostatic moment. For all values of the set of parameters, excluding some hypersurfaces, it is proved that the degenerate local and semi-local singularities of the system are of the parabolic and cuspidal type, respectively. Thus these singularities are structurally stable under small perturbations of the system in the class of integrable systems. Note that all these degenerate points lie in the preimage of the cusp point of the parametric bifurcation curve and satisfy the criterion of A.~Bolsinov, L.~Guglielmi and E.~Kudryavtseva.
Bibliogarphy: 28 items.
\end{abstract}


\textbf{1. Топологические инварианты интегрируемых систем.} Изучение гамильтоновых систем путем описания топологических свойств их поверхностей уровня энергии $H$ в фазовом пространстве было начато в известной работе С.~Смейла [1]. Как оказалось, вопрос о классе гомеоморфности таких изоэнергетических поверхностей при разных значениях энергии и других параметров системы (например, значения ее интеграла площадей $f_2$) тесно связан с результатами и подходами из теории Морса.

Если система интегрируема, то есть обладает ``достаточно большим'' числом скрытых симметрий, то изоэнергетическая поверхность сама может быть расслоена на совместные уровни дополнительных интегралов системы. На фазовом пространстве системы (являющемся симплектическим многообразием) возникает структура слоения Лиувилля, подробнее см. монографию [2].

Для таких систем возможен более тонкий топологический анализ. Слоения Лиувилля многих известных систем, возникающих в приложениях, (например, волчок Ковалевской, системы Жуковского и Горячева--Чаплыгина) были изучены М.П.~Харламовым. Для них были построены бифуркационные диаграммы отображения момента, определены классы гомеоморфности слоев и типы перестроек регулярных торов, см. работу [3].

Общая теория топологической классификации слоений Лиувилля интегрируемых систем основана на аналоге теории Морса, построенном для интегрируемых систем в работах А.Т.~Фоменко [4--5]. В дальнейших работах А.Т.Фоменко, его учеников и соавторов [6--8] были построены классифицирующие инварианты таких слоений относительно различных эквивалентностей.

Развитые подходы удалось применить к широкому классу интегрируемых систем из приложений. Для многих из них удалось вычислить классифицирующие инварианты Фоменко и Фоменко--Цишанга и определить типы невырожденных особенностей. Например, работы [9--14] посвящены изучению ряда систем динамики твердого тела, а недавние работы [15--17] --- аналогам таких систем на алгебрвх Ли, отличным от $e(3)$.

Интегрируемые системы разной природы могут оказаться эквивалентными. Например, их слоения Лиувилля в некоторых зонах энергии могут оказаться послойно гомеоморфными. Это означает одинаковое устройство замыканий решений этих систем в соответствующих зонах энергии. Такой эффект, в частности, позволяет ``моделировать'' поведение более сложной системы при помощи более простой и наглядной. В самые последние годы были получены интересные результаты по моделированию интегрируемых систем из механики при помощи интегрируемых биллиардов на кусочно-плоских столах, см. работы [18--21].

\textbf{2. Классификация особенностей интегрируемых систем.} Как оказалось, почти все критические точки (точки, где не максимален ранг дифференциала отображения момента) изученных конкретных систем невырождены, то есть удовлетворяют условию Морса--Ботта, см. [2, 5].

Особенности, все критические точки которых невырождены, удалось классифицировать. Классификация \textit{локальных} особенностей (то есть ростков слоения Лиувилля в окрестности точки) с точностью до симплектоморфизма была получена Элиассоном [22]. Класcы послойной гомеоморфности \textit{полулокальных} особенностей (то есть слоений в малой окрестности особого слоя, целиком состоящей из слоев слоения), коранга 1 были описаны А.Т.~Фоменко [2; теорема 3.3, т.1]. Их также называют 3-атомами Фоменко. Позже классы послойной гомеоморфности невырожденных особенностей любого коранга были описаны Н.Т.~Зунгом [23] с помощью подхода почти прямых произведений. Так были названы произведения ``элементарных'' полулокальных особенностей систем с одной или двумя степенями свободы, которые затем, возможно, факторизованы по действию конечной группы с определенными свойствами, подробнее см. [2; гл. 9 т.1.].

В большинстве систем с двумя степенями свободы из механики и математической физики, для которых была изучена топология слоений Лиувилля, почти все точки бифуркационной диаграммы (образа критического множества) имели в прообразе особые слои с невырожденными критическими точками ранга 1. В то же время, прообраз конечного числа точек, обыкновенно, устроен сложнее: в нем или имеются точки ранга 0 (невырожденные при почти всех значениях параметров) или одномерные орбиты, состоящие из вырожденных критических точек ранга 1.

Вопрос классификации вырожденных особенностей интегрируемых систем является весьма непростым. Так в работе А.В.~Болсинова, Л.~Гуглиелми и Е.А.~Кудрявцевой [24] был рассмотрен простейший класс таких локальных особенностей, называемых параболическими (также см. работу [25]). Локальная бифуркационная диаграмма таких особенностей имеет в подходящих координатах вид полукубической параболы.

В работе [24] был предложен критерий параболичности локальной вырожденной особенности ранга 1, сформулированный в терминах первых интегралов системы. Этот критерий оказался весьма удобным при проверке параболичности конкретных интегрируемых систем из приложений.

Если критическое множество на компактном особым слое состоит ровно из одной окружности, являющейся орбитой параболической точки (отметим, что тогда и остальные точки этой окружности тоже будут параболическими), то соответствующая полулокальная особенность называется \textit{каспидальной}, и является структурно устойчивой (т.е. неустранима при малом интегрируемом возмущении системы, то есть набора из симплектической структуры и первых интегралов, см. [24, 26]).

Предложенный в [24] критерий параболичности точки можно принять за определение таких точек. На симплектическом многообразии $(M^4, \Omega)$ рассмотрим отображение момента $\mathfrak{F} = (H, F)$, т.е. пару вещественно-аналитических функций $H$ и $F$ из $M^4$ в $\mathbb{R}$, которые коммутируют относительно формы $\Omega$. Такой набор $(H, F, \omega)$ определяет локальное гамильтоново $ \mathbb{R}^2$-действие на $M^4$. На неособых симплектических 4-листах рассматриваемых в нашей работе систем оно будет глобальным.

Положим, что $dF(x) \ne 0 $, и ограничим $H$ на поверхность уровня $F$, содержащую $x$, т.е. $H_0: = H|_{\{F = F(x)\}}$. Это гладкое 3-подмногообразие в $ M^4$. Также положим, что $ \mathrm{rk} \, d \mathfrak{F} = 1 $ в точке $x$. Тем самым $x$ --- критическая точка $H_0$, и $\exists$ единственное $k \in \mathbb {R} $, т.ч. $ dH(x) = k dF (x) $.

\textbf{Определение.}
\textit{
Точка $x$ (и проходящая через нее $ \mathbb{R}^2$-орбита гамильтонова действия) называется \textbf{параболической}, если выполнены три следующих условия: \\
(i) $\quad$ второй дифференциал ограничения $d^2 H_0(x)$ имеет ранг 1 в ней; \\
(ii) $\quad$ существует вектор $v \in \mathrm{Ker}\, d^2 H_0(x)$ т.ч. $v^3 H_0\ne 0$ (обозначим через $v^3H_0$ третий дифференциал $H_0$ в направлении касательного вектора $v$ в точке $x$); \\
(iii) $\quad$ второй дифференциал $d^2 (H-kF)(x)$ имеет ранг 3, где $k \in \mathbb{R}$ т.ч. $dH(x)\,=\,k\,dF(x)$.
}

Подчеркнем, что определение невырожденных точек через условие Морса--Ботта и их последующее активное изучение было обусловлено тем, что именно такими являются почти все особые точки систем из приложений.

Вполне естественным является вопрос: удовлетворяет ли определению параболичности достаточно широкий класс вырожденных особых точек реальных интегрируемых систем? В настоящей работе мы покажем, что для одного известного многопараметрического семейства систем динамики твердого тела свойство параболичности действительно является свойством общего положения.


\textbf{3. Интегрируемый случай Жуковского.} Задача о движении волчка, то есть о вращении твердого тела вокруг своей неподвижной точки, может быть обобщена и дополнительно усложнена добавлением к телу постоянного гиростатического момента. Как оказалось, при таком обобщении волчков Эйлера, Лагранжа и Ковалевской [28] получаются системы (они перечислены, например, в монографии [2]), которые тоже интегрируемы. Они известны в литературе, например, как случай Жуковского, случай Лагранжа с гиростатом и случай Ковалевской--Яхьи [14].

Интегрируемое обобщение случая Эйлера, т.е. тела, закрепленного на шарнире в своем центре масс, было открыто Н.Е.~Жуковским, см. [27]. Мы рассмотрим случай осесимметричного тела. Параметрами такой системы являются значение интеграла площадей $f_2 = b$, отношение двух различных моментов инерции $A_1:A_2$ и значения проекций $\lambda_1, \lambda_2$ вектора гиростатического момента на ось симметрии и на перпендикулярную ей плоскость.

Систему Жуковского, как и многие другие системы динамики твердого тела, можно задать на двойственном пространстве к алгебре Ли $e(3)$, т.е. на $\mathbb{R}^6$ с координатами $(J_1, J_2, J_3, x_1, x_2, x_3)$ и скобкой Ли--Пуассона
\[\{J_i, J_j\} = \varepsilon_{ijk} J_k, \qquad \{J_i, x_j\} = \varepsilon_{ijk} x_k \qquad \{x_i, x_j\} = 0,\]
где $\varepsilon_{ijk}$ есть знак перестановки $(123) \to (ijk)$. Сама динамическая система задается уравнениями вида $\dot{y_i} = \{y_i, H\}$ для $y_i \in \{J_1, \dots, x_3\}$.

Рассмотренная выше скобка имеет следующие функции Казимира, являющиеся первыми интегралами системы: геометрический интеграл $f_1$ и интеграл площадей $f_2$
\[f_1 = x_1^2 + x_2^2 + x_3^2, \qquad f_2 = x_1 J_1 + x_2 J_2 + x_3 J_3.\]
Совместная поверхность уровня $M^4_{a, b} : \{f_1 = a, f_2 = b\} \subset \mathbb{R}^6$ является гладким симплектическим многообразием при $a >0$ и гомеоморфна $T^* S^2$. Без ограничения общности положим $a = 1$, чего можно добиться линейной заменой координат. Обозначим  многообразие $M^4_{1, b}$ через $M^4_b$.

Гамильтониан системы Жуковского, умноженный на $2$, обозначим $H$ и приведем ниже. Он зависит от набора физических параметров: соотношения моментов инерции относительно главных осей $A_1 : A_2 : A_3$ и вектор гиростатического момента $\vec{\lambda} = (\lambda_1, \lambda_2, \lambda_3) \in \mathbb{R}^3$.
\[H = \cfrac{(J_1 + \lambda_1)^2}{A_1} + \cfrac{(J_2 + \lambda_2)^2}{A_2} + \cfrac{(J_3 + \lambda_3)^2}{A_3} \, ;\]
При произвольных значениях параметров (из физического смысла задачи $A_i > 0$) этот гамильтониан интегрируем, причем интеграл совпадает с интегралом системы Эйлера $F$ и является квадратичным:
\[F = J_1^2 + J_2^2 + J_3^2.\]
Также для удобства записи вычислений введем величины $a_i, \alpha_i$ и $\mu_i$ для $i = 1, 2, 3$
\begin{center}$a_i = 1 / A_i, \qquad a_i = \alpha_i^3,\qquad \lambda_i = \mu_i^3.$
\end{center}
В случае, если все момента инерции различны, а компоненты вектора $\vec{\lambda}$ отличны от нуля, фазовая топология системы была подробно изучена в работах М.П.~Харламова [3] и А.А.~Ошемкова [9]. При этом бифуркационная диаграмма отображения момента на плоскости $Ohf$ значений интегралов $H, F$ содержит параметрическую кривую $(h(t), f(t))$ для $t \in \mathbb{R}$:
\[h(t) = t^2 \left( \,
\cfrac{a_1 \lambda_1^2}{(a_1 - t)^2} +
\cfrac{a_2 \lambda_2^2}{(a_2 - t)^2} +
\cfrac{a_3 \lambda_3^2}{(a_3 - t)^2}
\,\right)\]
\[k(t) = \cfrac{a_1^2 \lambda_1^2}{(a_1 - t)^2} +
\cfrac{a_2^2 \lambda_2^2}{(a_2 - t)^2} +
\cfrac{a_3^2 \lambda_3^2}{(a_3 - t)^2}.
\]
При $a_1 > a_2 > a_3$ полученная кривая имеет три непрерывных ветви, соответствующие промежуткам $(a_3, a_2), (a_2, a_1), (a_1, +\infty) \cup (- \infty, a_3)$. Последние разделяются значениями параметра $t \in \{a_1, a_2, a_3\}$. Каждому из двух конечных интервалов соответствует кусочно-гладкая ветвь кривой, имеющая ровно одну точку возврата. Вычисление значения параметра кривой, соответствующего этой точке, для произвольных моментов инерции и компонент вектора гиростатического момента представляет собой непростую задачу, и ранее не выполнялось.

Система Жуковского, как и система Эйлера, обладает замечательным свойством. Поскольку гамильтониан $H$ и интеграл $F$ не зависят от переменных $x_1, x_2, x_3$, то $M^4_{b}$ можно спроецировать на $\mathbb{R}^3(J_1, J_2, J_3)$ и получить следующее слоение с особенностями. Каждой тройке $\vec{J} = (J_1, J_2, J_3)$ соответствуют точки, у которых набор координат $(x_1, x_2, x_3)$ как точка $\vec{x}$ пространства $\mathbb{R}^3(x_1, x_2, x_3)$ лежит в пересечении сферы $f_1 = |\vec{x}|^2 = a = 1$ и плоскости $b = <\vec{x}, \vec{J}>$.

При этом у точек открытого шара $|\vec{J}| < |b|$ прообраз в $M^4_{b}$ пуст, у точек сферы $|\vec{J}| = |b|$ прообраз состоит из точки $\vec{x} = \pm \vec{J}$ с нужными знаками, а каждой точке из внешности замкнутого шара $|\vec{J}| > |b|$ соответствует окружность.

Напомним, что $F = |\vec{J}|^2 = f$. Тем самым в $M^4_b$ все слои слоения Лиувилля с условием $f < b^2$ будут пусты, а в случае $f = b^2$ слой будет гомеоморфен своей проекции --- связной компоненте пересечения уровней $H = h$ и $F = f$ в $\mathbb{R}^3(J_1, J_2, J_3)$. В случае $f > b^2$ имеем, что каждый слой гомеоморфен произведению окружности и связной компоненты уровня $H = h, F = f$. Это произведение является прямым: топология слоев слоения Лиувилля для системы Жуковского хорошо известна. В частности, слоение Лиувилля не содержит бутылок Клейна или 3-атомов Фоменко со звездочками, т.е. невырожденных особенностей с неориентируемой сепаратрисной диаграммой.

\textbf{4. Основные результаты.} Рассмотрим случай осесимметричного тела. Тогда положим $A_2 = A_3$, то есть $a_2 = a_3$. Повернув, при необходимости, вторую и третью ось, примем $\lambda_3 = 0$, т.е. осесимметричная система зависит от четырех параметров: $a_1 \ne a_2, \lambda_1 \ne 0, \lambda_2 \ne 0$.

С точки зрения слоения Лиувилля и бифуркационной диаграммы при этом происходит следующее. У параметрической кривой $h(t), f(t)$ исчезает одна ветвь, содержащая одну из двух (в не осесимметричном случае Жуковского) точек возврата системы и соответствующая интервалу $(a_3, a_2) \subset \mathbb{R}$ ее параметра $t$. Оставшаяся часть бифуркационной диаграммы никаких перестроек не испытывает.

В осесимметричном случае несложно найти значение параметра $t_0$, соответствующее сохранившейся точке возврата кривой и вычислить координаты ее образа $h_0 = h(t_0)$ и $f_0 = f(t_0)$ на плоскости $Ohf$:
\[t_0 = \alpha_1^2 \alpha_2^2\cfrac{(\alpha_2 \mu_1^2 + \alpha_1 \mu_2^2)}{
\alpha_1^2 \mu_1^2 + \alpha_2^2 \mu_2^2}\,,
\qquad
h_0 = \alpha_1^3 \alpha_2^3\cfrac{(\alpha_2 \mu_1^2 + \alpha_1 \mu_2^2)^3}{(\alpha_1^3 - \alpha_2^3)^2}\,,
\quad
f_0 = \cfrac{(\alpha_1^2 \mu_1^2 + \alpha_2^2 \mu_2^2)^3}{(\alpha_1^3 - \alpha_2^3)^2}\,.
\]

Более трудной оказывается задача вычисления координат $(J_1, J_2, J_3)$, которые соответствовали бы критической окружности (состоящей из вырожденных точек).

\textbf{Утверждение. } 
В случае $|b| < b_0$, где $b_0 = +\sqrt{f_0} = (\alpha_1^2 \mu_1^2 + \alpha_2^2 \mu_2^2)^{3/2} / |\alpha_1^3 - \alpha_2^3|$, прообраз точки возврата кривой $h(t), f(t)$ при отображении момента $(H, F)$ гомеоморфен двумерному тору. Критическое множество на нем гомеоморфно окружности, и все ее точки вырождены. Эта критическая окружность проецируется в точку $(J_{10}, J_{20}, J_{30})$ пространства $\mathbb{R}^3(J_1, J_2, J_3)$ с координатами
\[J_{10} = -\alpha_1 \mu_1 \cfrac{\alpha_1^2 \mu_1^2 + \alpha_2^2  \mu_2^2}{\alpha_1^3 - \alpha_2^3},
\quad J_{20} = \alpha_2 \mu_2 \cfrac{\alpha_1^2 \mu_1^2 + \alpha_2^2 \mu_2^2}{\alpha_1^3 - \alpha_2^3},
\quad J_{30} = 0.\]

Далее несложно выбрать тройку координат $(x_1, x_2, x_3)$, задающую вместе с тройкой $(J_{10}, J_{20}, J_{30})$ вырожденную критическую точку в $M^4_b$.

Теперь применим критерий параболичности критической точки ранга 1 интегрируемой гамильтоновой системы с 2 степенями свободы, предложенный А.В.~Болсиновым, Л.~Гуглиелми и Е.А.~Кудрявцевой в [24]. Сформулируем основную теорему нашей работы:

\textbf{Теорема. } \textit{Рассмотрим осесимметричную систему Жуковского ($A_1 = A_2 < A_3$, либо $A_1 < A_2 = A_3$) с вектором гиростатического момента $\vec{\lambda} = (\lambda_1, \lambda_2, \lambda_3)$, таким что $\lambda_1 \ne 0, \lambda_2^2 + \lambda_3^2 \ne 0$. При условиях $f_1 = 1$ и $|f_2| = |b| < b_0$ все вырожденные особые точки в $M^4_{1, b}$ являются параболическими и составляют одну окружность. Содержащая ее полулокальная особенность является каспидальной и структурно устойчива при малых возмущениях системы в классе интегрируемых систем. }


\textbf{5. Доказательства. } Найдем координаты $(J_{10}, J_{20}, J_{30})$ точки, являющейся проекцией вырожденной окружности из прообраза точки возврата кривой $h(t), f(t)$.

\textbf{Доказательство утверждения}  1. Рассмотрим задачу нахождения условного экстремума функции $H$ на множестве $F = f$. Записав соответствующие условия на производные $H - \lambda F$ по $J_1, J_2, J_3$, имеем для некоторого $\lambda$
\[J_1 = -\cfrac{\mu_1^3}{1 + A_1 \lambda}, \qquad J_2 = -\cfrac{\mu_2^3}{1 + A_2 \lambda}, \qquad J_3 = 0.\]
Основную сложность представляет нахождение $\lambda = \lambda_0$ при котором получаемые критические точки попадают на уровень $F = f_0$. Подставив полученные выше значения в уравнения $H$ и $F$, имеем
\begin{equation}\lambda^2 \left(\cfrac{\alpha_1^6 \mu_1^6}{(1 - \alpha_1^3 \lambda)^2} + \cfrac{
   \alpha_2^6 \mu_2^6}{(1 - \alpha_2^3 \lambda)^2} \right) \, = \, \cfrac{(\alpha_1^2 \mu_1^2 + \alpha_2^2 \mu_2^2)^3}{(\alpha_1^3 - \alpha_2^3)^2} \, = f_0,
\end{equation}
\begin{equation}\cfrac{\alpha_1^3 \mu_1^6}{(1 - \alpha_1^3 \lambda)^2}  +
  \cfrac{\alpha_2^3 \mu_2^6}{(1 - \alpha_2^3 \lambda)^2} \, = \,
\alpha_1^3 \alpha_2^3\cfrac{(\alpha_2 \mu_1^2 + \alpha_1 \mu_2^2)^3}{(\alpha_1^3 - \alpha_2^3)^2} \, = h_0.
\end{equation}
Рассмотрим линейную комбинацию $f_0 - h_0 (2 \lambda - 1 / \alpha_1^3)$ с коэффициентом, зависящим от $\lambda$. После небольших упрощений получим выражение вида
\begin{equation} \cfrac{\mu_2^6 (\alpha_2^3 - \alpha_1^3)}{\alpha_1^3 (- 1 + \alpha_2^3 \lambda)^2} =
f_0 - h_0 \left(2 \lambda - \cfrac{1}{\alpha_1^3}\right) - (\mu_1^6 + \mu_2^6).
\end{equation}
Теперь выразим из уравнения (2) дробь $\cfrac{\alpha_1^3 \mu_1^6}{(1 - \alpha_1^3 \lambda)^2}$, а из уравнения (3) --- дробь $\cfrac{1}{(- 1 + \alpha_2^3 \lambda)^2}$. Подставим их в уравнение (1). Перенеся все в левую часть и избавившись от знаменателя, имеем
\begin{equation} \lambda^2 ((\alpha_2^3 - \alpha_1^3) \mu_2^6 \alpha_1^3 h_0 + (\alpha_2^3 - \alpha_1^3) \alpha_2^3 \mu_2^6 \alpha_1^3 (f_0 - h_0 (2 \lambda - 1/ \alpha_1^3) \\ - (\mu_1^6 + \mu_2^6))) -
 f_0 (\alpha_2^3 - \alpha_1^3) \mu_2^6 \, = 0.
\end{equation}
Уравнение (3) после домножения на знаменатель и уравнение (4) задают два многочлена третьей степени по $\lambda$. Рассмотрим их линейную комбинацию с коэффициентами $\alpha_2^3$ и $(-\alpha_1^3 + \alpha_2^3) \mu_2^6$ соответственно. Ее  старший коэффициент будет нулевым. Подстановка значений $h_0$ и $f_0$, найденных ранее, и последующее символьное упрощение с использованием системы Wolfram Mathematica 12 приводит к равенству вида: квадрат линейного по $\lambda$ многочлена равен нулю. Отсюда имеем единственное значение $\lambda_0$, соответствующее искомой критической точке (несложно убедиться, что эта точка не является морсовской на сфере $|J|^2 = f_0$):
\[\lambda = \lambda_0 = \cfrac{\alpha_1^2 \mu_1^2 +
 \alpha_2^2 \mu_2^2}{\alpha_1^2 \alpha_2^2 (\alpha_2 \mu_1^2 + \alpha_1 \mu_2^2)}.\]
Отсюда легко находим $J_{10}$ и $J_{20}$. Утверждение доказано. $\square$

\textit{Замечание. }
Отметим, что если $A_1 < A_2$, то $\alpha_1 > \alpha_2$ и $J_1 < 0, J_2 > 0$. Если же $A_1 > A_2$, то $\alpha_1 < \alpha_2$ и $J_1 > 0, J_2 < 0$. Это наблюдение будет полезно далее, когда при введении локальных координат $(J_1, J_2, x_1, x_2)$ на $M^4_b$ потребуется выбирать знак перед радикалами.

Теперь найдем некоторую точку из $M^4_b$, которая проецируется на найденную точку $(J_{10}, J_{20}, 0)$, и затем проверим ее параболичность.

\textbf{Доказательство теоремы.} 1. Введем на $M^4_{1, b}$ локальные координаты $J_1, J_2, x_1, x_2$, т.е. выразим $x_3$ и $J_3$ как функции от них из уравнений $f_1 = 1, f_2 = b$. Якобиан $\cfrac{\partial(f_1, f_2)}{\partial (J_3, x_3)}$ равен $2 x_3^2$, т.е. для регулярности координат достаточно брать точку из окружности, в которой $x_3  = \sqrt{1 - x_1^2 - x_2^2} \ne 0$.

Рассмотрим на $M^4_{1, b}$ первый интеграл $\Phi = H - \alpha_2^3 F$ и поверхность его уровня $Q^3$, содержащую слой с вырожденной окружностью. В качестве локальных координат выберем $J_2, x_1, x_2$, т.е. выразим $J_1$ через них. Якобиан замены имеет вид $4 \, \alpha_1 \alpha_2^2 \, \mu_1 \mu_2^2 \, (- 1 + x_1^2 + x_2^2)$, т.е. отличен от нуля при $x_3 \ne 0$ и $\mu_1 \mu_2 \ne 0$.

2. Плоскость $b = <\vec{x}, \vec{J_0}>$ параллельна оси $Ox_3$ в $\mathbb{R}^3(x_1, x_2, x_3)$, поскольку $J_{30} = 0$. Выберем в качестве $x_1, x_2$ значения координат в точке основания перпендикуляра, опущенного на эту плоскость из начала координат. Их легко найти из вида $J_{10}$ и $J_{20}$. Более того, значение $x_3 = \sqrt{1 - x_1^2 - x_2^2}$ обязательно положительно в случае, если пересечение плоскости и сферы $f_1 = |\vec{x}|^2 =  1$ является окружностью.
\[x_{10} = -\cfrac{\alpha_1 \mu_1 (\alpha_1^3 - \alpha_2^3) b}{(\alpha_1^2 \mu_1^2 + \alpha_2^2 \mu_2^2)^2}, \qquad x_{20} = +\cfrac{\alpha_1 \mu_1 (\alpha_1^3 - \alpha_2^3) b}{(\alpha_1^2 \mu_1^2 + \alpha_2^2 \mu_2^2)^2}, \qquad x_{30} = +\sqrt{1 - x_{10}^2 - x_{20}^2}.\]

3. Выразив $x_3$ из уравнения $f_1 = 1$, и подставив в $f_2 = b$, имеем
\[x_3 = \sqrt{1 - x_1^2 - x_2^2}, \qquad J_3 = \cfrac{b - x_1 J_1 - x_2 J_2}{\sqrt{1 - x_1^2 - x_2^2}}.\]
Подставив эти выражения в формулу интеграла $\Phi$, имеем
\[\Phi = \alpha_1^3 (J_1 + \mu_1^3)^2 + \alpha_2^3 (-J_1^2 + 2 J_2 \mu_2^3 + \mu_2^6).\]
Значение $\Phi = \phi_0$ в вырожденных точках равно
\begin{center}
$\varphi_0 = \cfrac{\alpha_2^3 (3 \alpha_1^2 \alpha_2 \mu_1^2 \mu_2^4 + \alpha_2^3 \mu_2^6 +
   \alpha_1^3 (-\mu_1^6 + \mu_2^6))}{\alpha_1^3 - \alpha_2^3}$
   \end{center}
Получили многочлен второй степени по переменной $J_1$. Выражая $J_1$, получим две ветви, отличающиеся знаком перед радикалом. Подставив значение $J_{20}$, имеем
\begin{center}
$J_1 = \cfrac{-\alpha_1^3 \mu_1^3 \pm \alpha_1 \alpha_2^2 |\mu_1| \mu_2^2}{\alpha_1^3 - \alpha_2^3}.$
\end{center}
Требуется выбрать знак, при котором полученный $J_1$ равен $J_{10}$. Без ограничения общности примем $\mu_1 >0$ (иначе поменяем направление оси). Тогда надо выбрать знак ``минус''. Получаем
\[J_1(J_2, x_1, x_2) \,\, = \,\,
- \cfrac{\alpha_1^3 \mu_1^3 + \sqrt{\alpha_2^3 \mu_2^3} \, \sqrt{2 J_2 (\alpha_2^3 - \alpha_1^3) + 3 \alpha_1^2 \alpha_2 \mu_1^2 \mu_2 +
     2 \alpha_2^3 \mu_2^3}}{(\alpha_1^3 - \alpha_2^3)}.\]

4. Рассмотрим найденный выше интеграл $\Phi(J_1, J_2, x_1, x_2)$, а также функцию $F(J_1, J_2, J_3(J_1, J_2, x_1, x_2))$ и вычислим их дифференциалы по переменным $(J_1, J_2, x_1, x_2)$. Выполнив символьные вычисления в системе Wolfram Mathematica 12, подставим в эти формулы значения координат $J_{10}, J_{20}, x_{10}, x_{20}$ вырожденной критической точки. Получим, что дифференциал следующей линейной комбинации обращается в ноль:
\begin{center}
$(\alpha_1^2 \mu_1^2 + \alpha_2^2 \mu_2^2) d \Phi  - \alpha_2^2 \mu_2^2 (\alpha_1^3 - \alpha_2^3) d F = (0, 0, 0, 0).$
\end{center}

Рассмотрим второй дифференциал этой линейной комбинации. Аналогично, используя систему Wolfram Mathematica 12, вычислим символьно соответствующую матрицу, и затем подставим координаты изучаемой точки.

Две строки полученной матрицы, соответствующие переменным $x_1$ и $x_2$, линейно зависимы. Тем самым, для выполнения условия (iii) определения параболической точки остается проверить невырожденность подматрицы $3 \times 3$, соответствующей переменным $J_1, J_2, x_1$.

Ее определитель есть произведение $8 \,\, \alpha_1^4 \alpha_2^4 \,\, (\alpha_1^3 - \alpha_2^3) \,\, \mu_1^4 \,\, \mu_2^4 \,\, (\alpha_1^2 \mu_1^2 + \alpha_2^2 \mu_2^2)^5$ и квадрата следующего выражения от $b$:
\begin{center} $-2 \alpha_1^3 \alpha_2^3 b^2 + \alpha_1^6 (b^2 - \mu_1^6) - 3 \alpha_1^4 \alpha_2^2 \mu_1^4 \mu_2^2 - 3 \alpha_1^2 \alpha_2^4 \mu_1^2 \mu_2^4 +  \alpha_2^6 (b^2 - \mu_2^6).$
 \end{center}
Последнее отлично от нуля, если $|b| \ne f_0$ --- значению $F$ на вырожденной орбите. Тем самым, при $|b| < \sqrt{f}$ ранг матрицы второго дифференциала равен трем, т.е. это условие (iii) параболичности выполнено.

5. Остается проверить условия (i) и (ii). Сделаем это для $F$ как функции от локальных координат $J_2, x_1, x_2$ на совместной поверхности уровня функций $f_1, f_2, \Phi$ вблизи изучаемой точки. Полученная матрица второго дифференциала в данной точке имеет нулевыми первую строку и столбец (дифференцирование по $J_2$), а две другие строки линейно зависимы и не являются ненулевыми. Тем самым, ранг второго дифференциала ограничения равен 1, и условие (i) выполнено.

Для проверки условия (ii) заметим, что вектор $(1, 0, 0)$ лежит в ядре второго дифференциала. Подставим его в третий дифференциал этой функции $F$, и после упрощений получим $\cfrac{6 (\alpha_1^3 - \alpha_2^3)}{\alpha_1^2 \alpha_2 \mu_1^2 \mu_2}$. В случае $\alpha_1 \ne \alpha_2$ и $\mu_1 \mu_2 \ne 0$ он корректно определен и отличен от нуля, т.е. условие (ii) выполнено, и данная точка является параболической (вместе со всей ее орбитой). Теорема 1 доказана. $\square$

\end{fulltext}

\end{document}